\documentclass{amsart}

\usepackage{amsmath,amsthm,amssymb,amscd,accents,enumitem,mathtools,mathrsfs}
\usepackage{amsfonts}
\usepackage{rotating}
\usepackage{euscript}
\usepackage{pst-node}
\usepackage{epsfig,verbatim}
\usepackage{tikz-cd}
\def\be{\begin{equation}}
\def\ee{\end{equation}}

\def\C{{\mathbb C}} 
\def\f{\EuScript}
\def\oo{\mathscr}
 
\def\P{{\mathbb P}}

\def\Q{{\mathbb Q}}

\def\phi{{\varphi}}

\def\tt{\widetilde}
\def\deg{{\rm deg\,}} 

\def\Aut{{\rm Aut}}
\def\Mon{{\rm Mon}}

\def\GCD{{\rm GCD }}

\def\bp{\begin{proposition}}
\def\ep{\end{proposition}}

\def\bt{\begin{theorem}}
\def\et{\end{theorem}}
\def\br{\begin{remark}}
\def\er{\end{remark}}
\def\be{\begin{equation}}
\def\bee{\begin{equation*}}
\def\l{\label}

\def\ee{\end{equation}}
\def\eee{\end{equation*}}
\def\bl{\begin{lemma}}
\def\el{\end{lemma}}
\def\bc{\begin{corollary}}
\def\ec{\end{corollary}}
\def\pr{\noindent{\it Proof. }}

\def\bd{\begin{definition}}
\def\ed{\end{definition}}
\def\t{\widetilde}
\def\hat{\widehat}

\mathtoolsset{showonlyrefs}

\newtheorem{theorem}{Theorem}[section]
\newtheorem{lemma}[theorem]{Lemma}
\newtheorem{corollary}[theorem]{Corollary}
\newtheorem{proposition}[theorem]{Proposition}
\theoremstyle{definition}
\newtheorem{remark}[theorem]{Remark}

\begin{document}
\title[Lower bounds for genera of fiber products]{Lower bounds for genera of fiber products}

\begin{abstract}
We give lower bounds for genera of components of fiber products of holomorphic maps between compact Riemann surfaces, 
extending  results on genera of components of algebraic curves of the form 
$A(x)-B(y)=0,$ where $A$ and $B$ are rational functions. 

\end{abstract} 
\author{Fedor Pakovich}
\thanks{
This research was supported by ISF Grant No. 1092/22}
\address{Department of Mathematics, Ben Gurion University of the Negev, Israel}
\email{
pakovich@math.bgu.ac.il}

\maketitle

\section{Introduction}
In this paper, we extend results of the recent papers \cite{cur}, \cite{tame} concerning lower bounds for genera of components of algebraic curves of the form \be \l{cur} E_{A,B}:\, A(x)-B(y)=0,\ee where $A$ and $B$ are rational functions with complex coefficients, to the case of fiber products of arbitrary holomorphic maps between compact Riemann surfaces.  Not less  importantly, we simplify the approach used in \cite{cur}, \cite{tame} directly relating the  problem to the group action on Riemann surfaces and the Hurwitz automorphisms theorem.  Here and below, we always assume that considered functions and maps are non-constant.

The problem of describing rational functions $A$, $B$ 
such that the algebraic curve \eqref{cur} 
 has a factor of genus zero or  one, to which we refer below as ``the low genus problem'', naturally arises in several different branches of mathematics.

First, since \eqref{cur} has a factor of genus zero if and only if there exist rational functions $X,Y$ satisfying \be \l{2} A\circ X=B\circ Y,\ee the low genus problem is central in the theory of functional decompositions of rational functions. 
In the  polynomial case, this theory was developed  by Ritt 
 (see \cite{r1}, \cite{sch}). The general case, however, is much less understood and 
known results are mostly concentrated on a study of either decompositions of special types of functions or functional equations of a special form (see e.g.  \cite{az}, \cite{bog}, \cite{ere}, \cite{fz}, \cite{mp}, \cite{ngw}, \cite{lau}, \cite{semi}, \cite{rev}, \cite{ritt}). Notice also that  by the  Picard theorem  any algebraic curve that can be
parametrized by functions meromorphic on $\C$ has genus zero or  one. 
Thus, the functional equation \eqref{2}, where $X$, $Y$ are allowed to be entire or meromorphic functions on $\C$, often studied in the context of Nevanlinna theory (see e.g. \cite{vd0}, \cite{fuj}, \cite{vd1}, \cite{vd3}, \cite{yang}), 
is also related to the low genus problem 
(see e. g. \cite{bn}, \cite{vd2}, \cite{ent}, \cite{alg}).

Second, algebraic curves \eqref{cur} with factors of genus zero or one have special Diophantine properties. Indeed, by the  Siegel theorem, if an irreducible algebraic curve $C$ with rational coefficients has infinitely many integer points, then $C$ is of genus zero with at most two points at infinity. More generally, by the Faltings theorem,  if $C$ has infinitely many rational points, then  $g(C)\leq 1.$ Consequently,  since many interesting Diophantine equations have the form $ A(x)=B(y)$, where $A,B$ are rational functions over $\Q$, the low genus problem is of importance in the number theory (see e.g. \cite{az2}, \cite{bf}, \cite{bilu}, \cite{dk}, \cite{f1}, \cite{kr}, \cite{kt}, \cite{sto}).    
The most notable result in this direction is the complete classification of polynomial  curves $E_{A,B}$ having a factor of genus zero  with at most two  points at infinity
 obtained in the
paper of Bilu and Tichy \cite{bilu}, which continued the line of researches started by Fried (see \cite{f1}, \cite{f2}, \cite{f3}). 

Third, the low genus problem naturally arises in the new emerging field of arithmetic dynamics.  
For example, the problem of describing rational functions $A$ and $B$ such that all curves 
\be \l{cu1} A^{\circ n}(x)-B(y)=0, \ \ \ n\geq 1,\ee
 have a factor of genus zero or one 
is a geometric counterpart of the following problem of the arithmetic nature (see \cite{jones}, \cite{hyde}, \cite{aol}):   for which rational functions $A$ and $B$ 
defined over a number field $K$ there exists a $K$-orbit of $A$ containing infinitely many points from the value set $B(\P^1(K))$  ? More generally, the problem of describing pairs of rational functions  $A$ and $B$ such that all curves  
\be \l{cu2} A^{\circ n}(x)-B^{\circ m}(y)=0, \ \ \ n,m\geq 1,\ee
have a factor of genus zero or one is a geometric counterpart of the problem of describing pairs of rational functions $A$ and $B$ having orbits
with infinite intersection (see \cite{gtz}, \cite{gtz2}, \cite{tame}). 

 Finally, notice that the low genus problem is related to the study of amenable semigroups of rational functions  under the operation of composition,
 since for such a semigroup  $S$ the amenability condition implies that for all $A,B\in S$ all curves \eqref{cu2} have a factor of genus zero    
 (see \cite{amen}).

In case  the curve $E_{A,B}$ is irreducible,
 the standard approach to the low genus problem   initiated by Fried (\cite{f1}, \cite{f3}) relies on the use of   
 an explicit formula for genus of $E_{A,B}$ in terms of the ramifications of $A$ and $B$ (see Section \ref{fib} below).  However, the direct analysis of this formula is quite  difficult, and obtaining   a full classification of curves $E_{A,B}$ of genus zero or one on this way 
seems to be hardly possible. In addition, such an analysis 
results  only in possible {\it patterns} of ramifications of $A$ and $B$. However,  rational functions with such  patterns may not  exist.  It is known that any ``polynomial'' pattern is realizable by some polynomial (\cite{thom}), but 
already for ``Laurent polynomial'' patterns there exists a number of exceptions (\cite{hu}).  
In general, the problem of existence of a rational function with a given ramification pattern, the so-called Hurwitz problem, is still widely open (see e. g. the recent papers \cite{hu00}, \cite{hu0}, \cite{hu1}, \cite{hu2} and the bibliography therein). 

A general lower bound for the genus of $E_{A,B}$ was obtained in the paper \cite{cur}. To formulate it explicitly, 
let us recall that for a  holomorphic map between compact Riemann surfaces $P:\, \f R\rightarrow \f C$  its {\it normalization} is defined as a compact Riemann surface  $\f N_P$ together with a holomorphic Galois covering  of the 
lowest possible degree $\t P:\f N_P\rightarrow \f C$ such that
 $\t P=P\circ H$ for some  holomorphic map $H:\,\f N_P\rightarrow \f R$. From the algebraic point of view, the passage from $P$ to $\t P$ 
corresponds to the passage from the field extension $ \f M(\f R)/P^*\f M(\f C)$ to its Galois closure. 
In these terms, the main result of \cite{cur}  may be formulated as follows: if  
 $A$ and $B$ are rational functions of degree $n$ and $m$  correspondingly such that 
 $E_{A,B}$ is irreducible and $g(\f N_A)>1$, then  
\be \l{ma} 
g( E_{A,B}) >\frac{m-84n+168}{84}.
\ee   Thus, for fixed $A$ the genus of $E_{A,B}$ grows linearly with respect to the degree of $B$, unless $A$ satisfies the condition $g(\f N_A)\leq 1.$ In particular, 
 $E_{A,B}$ has genus greater than one whenever $m\geq 84 (n-1)$.  
What is important is that the condition $g(\f N_A)\leq 1$ is quite restrictive. Specifically, up to the change 
$A\rightarrow \alpha\circ A\circ \beta,$ where $\alpha$ and $\beta$ are M\"obius transformations, the list of rational functions $A$ with $g(\f N_A )=0$ consists of the series 
$$z^n, \ n\geq 1, \ \ \ \ T_n, \ n\geq 2,  \ \ \ \ \frac{1}{2}\left(z^n+\frac{1}{z^n}\right), \ n\geq 1,$$ and a finite number of functions, which can be calculated explicitly (see \cite{gen}). On the other hand, functions with $g(\f N_A)=1$ admit a simple geometric description. In particular, the simplest examples of such functions are  Latt\`es maps  (see \cite{gen}).

In case  the curve $E_{A,B}$ is reducible, the above mentioned genus formula 
cannot be used for studying the low genus problem.  On the other hand,  
 the problem of reducibility of $E_{A,B}$, the so-called Davenport-Lewis-Schinzel problem, is very subtle and difficult 
(see  \cite{CCN99}, \cite{DLS}, \cite{DS}, \cite{f2},   \cite{Fried5}, \cite{FG},  \cite{nef2},  \cite{ssch}). Thus,  universal bounds for genera of  {\it components} of $E_{A,B}$ 
are especially interesting. However, it is easy to see that it is not possible to provide such bounds for {\it all} components of $E_{A,B}$,  since  for  arbitrary rational functions $A$ and $S$,  setting $B=A\circ S$ we obtain a  curve $E_{A,B}$ with an irreducible component of genus zero $x-S(y)=0.$ Nevertheless, it was shown in \cite{tame}  
that excluding from the consideration components of the above form and changing the condition $g(\f N_A)>1 $ to a stronger condition makes the problem solvable.

To formulate the  result of \cite{tame} explicitly, let us introduce the following definition. We 
say that a rational function $A$ is {\it tame} if the algebraic curve 
$$A(x)-A(y)=0$$ has no factors of genus zero or one distinct from the diagonal $x-y=0.$  Otherwise, we say that $A$ is {\it wild}. It can be shown that for every tame rational function $S$ the inequality  $g(\f N_A)> 1$ holds (\cite{tame}).  
Thus, the tameness condition is a strengthening of the  condition $g(\f N_A)> 1$.  Notice that a generic rational function of degree at least four is tame (\cite{alg}), but a comprehensive classification of  wild rational functions is not known (for some partial results see  \cite{ann}, \cite{az}, \cite{alg}, \cite{me},  \cite{r4},    \cite{seg}). 
In this notation, the result of \cite{tame} 
 can be formulated as follows: if $A$ is a tame rational function of degree $n$  
and $B$ a rational function 
 of degree $m$, then for any irreducible component $C$ of the curve $E_{A,B}$ the inequality    
\be \l{nera} 
g( C) \geq \frac{m/n!-84n+168}{84} \ee holds,  
unless  $B=A\circ S$ for some rational function $S$, and $C$ is the graph $x-S(y)=0.$

The algebraic curve $E_{A,B}$ can be interpreted as the {\it fiber product} of rational functions $A$ and $B,$ and in this paper we generalize results of \cite{cur}, \cite{tame} to the fiber products of arbitrary holomorphic maps between compact Riemann surfaces (see Section \ref{fib} for precise definitions). In practical terms, we consider commutative diagrams
\be \l{dii}
\begin{CD}
\f E @>U>> \f T\\
@VV {V} V @VV W V\\ 
\f R @>P >> \f C
\end{CD}
\ee
consisting of holomorphic maps between compact Riemann surfaces subject to the condition that  the maps $V$ and $U$ have no {\it non-trivial common compositional  right factor} in the following sense: 
the equalities 
$$U= \tt U\circ  T, \ \ \ V= \tt V\circ T,$$ where $$T:\, {\f E} \rightarrow \t{\f E} , \ \ \ \tt V:\, \t{\f E} \rightarrow \f R, \ \ \ \tt U:\,  \t{\f E}  \rightarrow \f T$$ are holomorphic maps between compact Riemann surfaces, imply that $\deg T=1.$  
For brevity, we will call such diagrams {\it reduced}. Notice that for a reduced diagram, the inequalities $\deg W\geq \deg V$ and $\deg P\geq \deg U$ hold.

 Our first main result  is the following generalization of the bound \eqref{ma} 
to the bound for the genus of the fiber product of  
holomorphic maps between compact Riemann surfaces, in case this product consists of a unique component.

\bt \l{t1} Let $P:\f R\rightarrow \f C$  and $W:\f T\rightarrow \f C$ be holomorphic maps between compact Riemann surfaces such that the fiber product $(\f R,P)\times_{\f C} (\f T,W)$ consists of a unique component $\f E$ and $g(\f N_W)> 1$. Then 
\be \l{svk4} g(\f E) \geq (g(\f R)-1)(\deg W-1)+1+\frac{\deg P}{84}.\ee
\et 

Notice that since $g(\f R)\geq 0$, inequality \eqref{svk4} implies the bound   
\be \l{svk0} g(\f E) \geq  \frac{\deg P-84\,\deg W+168}{84},\ee which depends only on the degrees of $P$ and $W$.  
Notice also that the condition $g(\f N_W)>1$ is obviously satisfied whenever $g(\f T)>1$. 

Let us mention that in general the bound provided by  Theo\-rem \ref{t1} is the best possible. Namely, the equality in \eqref{svk4} is attained for the fiber products of the quotient  maps associated with the action of full automorphism groups 
of the Hurwitz surfaces (see Section \ref{exam}). 

Finally, we remark that Theo\-rem \ref{t1} is not true if $g(\f N_W)\leq 1$ (see \cite{cur}). The simplest examples are  obtained from the commutative diagram 
\be \l{dur} 
\begin{CD}
\C\P^1 @>z^rR(z^n) >> \C\P^1\\
@VV {z^n} V @VV z^n V\\ 
\C\P^1 @>z^rR^n(z) >> \C\P^1,
\end{CD}
\ee
where $R$ is an arbitrary rational function and $r,n$ are integer positive numbers. Indeed, since the curve $E_{A,B}$ is irreducible whenever the degrees of $A$ and $B$ are coprime, choosing $R$ and $r,n$ appropriately we obtain an irreducible curve of genus 
zero $x^n-y^rR^n(y)=0$ such that for fixed $A=z^n$ the degree of $B=z^rR^n(z)$ can be arbitrary large.  

Our second main result is the following bound applicable to fiber products with several components.

\bt \l{t2}  Let 
\be 
\begin{CD} \l{dubi} 
\f E @>U>> \f T\\
@VV {V} V @VV W V\\ 
\f R @>P >> \f C
\end{CD}
\ee
be a reduced commutative diagram of holomorphic maps between compact Riemann surfaces   
such that $\deg V>1 $ and  the fiber product of $(\f T,W)$ with itself $\deg V$ times contains no components of genus zero or one that do not belong to the big diagonal of $\f T^{\deg V}$.
Then 
\be g(\f E) \geq (g(\f R)-1)(\deg V-1)+1\\+\frac{\deg P}{\deg W(\deg W-1)\dots (\deg W-\deg V+1)}.\ee  
\et 

Notice that the condition  $g(\f N_W)>1$ is equivalent to the condition that the fiber product of $(\f T,W)$ with itself $\deg W$ times contains no components of genus zero or one that do not belong to the big diagonal of $\f T^{\deg W}$ (see Section \ref{sec}). Hence,  in case \be \l{inc} \deg V=\deg W,\ee the assumption about $W$ from  Theorem \ref{t2} is equivalent to the assumption about $W$  from Theorem \ref{t1}. Moreover, if  \eqref{inc} holds, then the fiber product 
$(\f R,P)\times_{\f C} (\f T,W)$ consists of a unique component $\f E$ (see Section \ref{fib}), and 
thus the assumptions of both theorems coincide. Nevertheless, the bound provided by Theorem \ref{t2} 
is weaker than the bound provided by Theorem \ref{t2}. 

Applying  Theorem \ref{t2} to rational functions, we obtain the 
following statement.

\bt \l{rat} 
Let $A$ and $B$ be rational functions of degree $n$  and $m$ correspondingly, and $C:F(x,y)=0$ an irreducible component of the curve $E_{A,B}$ such that $k=\deg_xF>1$. Then 
\be \l{svk10} g(C) > 2-k+\frac{m}{n(n-1)\dots (n-k+1)},\ee
unless the algebraic curve in $(\C\P^1)^k$ defined by the equation 
\be \l{res} A(x_1)=A(x_2)=\dots =A(x_k)\ee has 
a component of genus zero or one that does not belong to the big diagonal of $(\C\P^1)^k$.
\et

Finally, we prove the following result slightly improving the result of \cite{tame}. 

\bt \l{ratt} 
Let $A$ and $B$ be rational functions of degree $n$  
and $m$ correspondingly, and $C$ an irreducible component  of the curve $E_{A,B}$. Assume that $A$ is tame. 
Then  
\be \l{svk11}
g(C) > 2-n+\frac{m}{n!}, \ee  
unless  $B=A\circ S$ for some rational function $S$, and $C$ is the graph $x-S(y)=0.$
\et

In brief, our proof of Theorem \ref{t1} goes as follows. First,  we establish a lower bound for 
 the Euler characteristic $\chi(\f N_V)$ of the normalization  of 
a holomorphic map between compact Riemann surfaces $V:\, \f E\rightarrow \f R$ in terms of $\chi(\f E)$ and $\chi(\f R)$,  and  $\deg V$ and $\deg \tt V$ (Section \ref{ces}). Using the Hurwitz automorphisms theorem, we also obtain 
an upper bound for $\chi(\f N_V)$ in  case $g(\f N_V)>1$.  
Then, we show that diagram  \eqref{dubi} can be lifted to a diagram of holomorphic maps between compact Riemann surfaces   
\be \l{medved} 
\begin{CD}
\f N_V @>L >> \f N_W \\
@VV Q V @VV F V\\ 
\f E @>U>> \f T\\
@VV {V} V @VV W V\\ 
\f R @>P >> \f C 
\end{CD}
\ee
(Section \ref{lii}).  
 Finally, we 
apply the Riemann-Hurwitz formula to the map \linebreak $L:\f N_V\rightarrow \f N_W$ and use the above mentioned bounds for 
$\chi(\f N_V)$ and $\chi(\f N_W)$. The proof of Theorem \ref{t2} is similar with the exception that in  diagram \eqref{medved} instead of 
$\f N_W$  appears  some irreducible component of the  fiber product  
of $W:\f T\rightarrow \f C$ with itself $\deg V$ times, and a rougher upper bound for the Euler characteristic of this component is used. 

In concluding this introduction, let us mention that an {\it upper} bound for the genera of components of fiber products follows from the  classical Castelnuovo-Severi inequality 
$$g(\f E)\leq g(\f R)\deg V+g(\f T)\deg U +(\deg V-1)(\deg U-1)$$  
for the genus 
of a compact Riemann surface $\f E$ such that there exist holomorphic maps 
$V:\f E\rightarrow \f R$ and $U:\f E\rightarrow \f T$ having no non-trivial common compositional  right 
factor (see \cite{acc1}, \cite{acc2}, \cite{kani}). Indeed, considering for a component of the fiber product of maps $P:\f R\rightarrow \f C$  and $W:\f T\rightarrow \f C$ the corresponding reduced diagram \eqref{dubi}, and taking into account that 
 $\deg W\geq \deg V$ and $\deg P\geq \deg U$, we conclude that 
$$g(\f E)\leq g(\f R)\deg W+g(\f T)\deg P +(\deg W-1)(\deg P-1).$$

The paper is organized as follows. In the second section, we recall some basic facts about fiber products and normalizations. 
In the third section, following the approach described above, we prove Theorems \ref{t1} - \ref{ratt}. We also show that the bound provided by Theorem  \ref{t1} is sharp.

\section{Fiber products and normalizations} 

\subsection{\l{fib} Fiber products}
Let $P:\, \f R\rightarrow \f C$ and $W:\, \f T \rightarrow \f C$ be holomorphic maps between compact Riemann surfaces. 
The collection
\be \l{nota} (\f R,P)\times_{\f C} (\f T,W)=\bigcup\limits_{j=1}^{n(P,W)}\{\f E_j,V_j,U_j\},\ee 
where $n(P,W)$ is an integer positive number and $\f E_j,$ $1\leq j \leq n(P,W),$ are compact Riemann surfaces provided with holomorphic maps
$$V_j:\, \f E_j\rightarrow \f R, \ \ \ U_j:\, \f E_j\rightarrow \f T, \ \ \ 1\leq j \leq n(P,W),$$
is called the {\it fiber product} of  $P$ and $W$ if $$ P\circ V_j=W\circ U_j, \ \ \ 1\leq j \leq n(P,W),$$ 
and for any holomorphic maps $V:\, \t{\f E}\rightarrow \f R,$  $U:\, \t{\f E}\rightarrow \f T$
between compact Riemann surfaces satisfying 
\be \l{pes}  P\circ V=W\circ U\ee there exist a uniquely defined  index $j$, $1\leq j \leq n(P,W)$, and 
a holomorphic map $T:\, \t{\f E}\rightarrow \f E_j$ such that
\be \l{ae} V= V_j\circ  T, \ \ \ U= U_j\circ T.\ee  

The fiber product  is defined in a unique way up to natural isomorphisms and 
can be described by the following algebro-geometric construction. Let us consider the algebraic variety 
\be \l{ccuurr} \f L=\{(x,y)\in \f R\times \f T \, \vert \,  P(x)=W(y)\}.\ee
Let us denote by $\f L_j,$ $1\leq j \leq n(P,W)$,  irreducible components of $\f L$, by 
$\f E_j$, \linebreak  $1\leq j \leq n(P,W)$, their desingularizations, 
 and by $$\pi_j: \f E_j\rightarrow \f L_j, \ \ \ 1\leq j \leq n(P,W),$$ the desingularization maps.
Then the compositions  $$x\circ \pi_j: \f E_j\rightarrow \f R, \ \ \ y\circ \pi_j: \f E_j\rightarrow \f T, \ \ \ 1\leq j \leq n(P,W),$$ 
extend to holomorphic maps
$$V_j:\, \f E_j\rightarrow \f R, \ \ \ U_j:\, \f E_j\rightarrow \f T, \ \ \ 1\leq j \leq n(P,W),$$
and the collection $\bigcup\limits_{j=1}^{n(P,W)}\{\f E_j,V_j,U_j\}$ is the fiber product of $P$ and $W$. 
Abusing notation we call the Riemann  surfaces $\f E_j$,  $1\leq j \leq n(P,W),$ irreducible components of the fiber product of $P$ and $W$.

It follows from the definition  that for every $j,$ $1\leq j \leq n(P,W),$ the functions $V_j,U_j$  have no {\it non-trivial common compositional  right factor} in the following sense: 
the equalities 
$$ V_j= \tt V\circ  T, \ \ \ U_j= \tt U\circ T,$$ where $$T:\, {\f E}_j \rightarrow \t{\f E} , \ \ \ \tt V:\, \t{\f E} \rightarrow \f R, \ \ \ \tt U:\,  \t{\f E}  \rightarrow \f T$$ are holomorphic maps between compact Riemann surfaces, imply that $\deg T=1.$  
Denoting by $\f M(\f R)$ the field of meromorphic functions on a Riemann surface $\f R$, we can restate  this  condition as the equality
$$ V_j^*\f M(\f R)\cdot U_j^*\f M(\f T)=\f M(\f E_j),$$ meaning that the field $\f M(\f E_j)$ is the compositum  of its subfields $V_j^*\f M(\f R)$ and $U_j^*\f M(\f T).$  
In the other direction, if $U$ and $V$ satisfy \eqref{pes} and have no non-trivial common compositional  right factor, then   
 $$V=V_j\circ T, \ \ \ \ U=U_j\circ T$$ for some $j$, $1\leq j \leq n(P,W),$ and an  isomorphism $T:\, \f E_j\rightarrow \f E_j.$  

Notice that since $V_i,U_i$, $1\leq j \leq n(P,W),$  parametrize components of   
\eqref{ccuurr}, the equalities    
\be \l{ii} \sum_{j}\deg V_j=  \deg W,\ \ \ \ \sum_{j}\deg U_j= \deg P\ee
hold. In particular, if $(\f R,P)\times_{\f C} (\f T,W)$ consists of a unique component $\{\f E,V,U\},$ then \be \l{vv} \deg V=\deg W,\ \ \ \ \deg U=\deg P.\ee Vice versa,   if holomorphic maps $U$ and $V$ satisfy \eqref{pes}  and \eqref{vv}, and have no non-trivial common compositional  right factor,  then $(\f R,P)\times_{\f C} (\f T,W)$ consists of a unique component.

If $R:\, \f E\rightarrow \f C$ is a holomorphic map between compact Riemann surfaces, then by the Riemann–Hurwitz formula  
\be \l{hur} \chi(\f E)=\chi(\f C)\deg R-\sum_{z\in \f E}(e_R(z)-1),\ee where $e_R(z)$ denotes the local multiplicity of $R$ at the point $z$. In case the fiber product 
$(\f R,P)\times_{\f C} (\f T,W)$ consists of a unique component $\{\f E,V,U\},$  one can calculate $\chi(\f E)$ applying \eqref{hur} to the map 
\be \l{dq} R=P\circ V=W\circ U\ee as follows (see \cite{f3} or \cite{lau}, Section 2). 

Setting $$R_j=P\circ V_j=W\circ U_j, \ \ \ 1\leq j\leq n(P,W),$$ let us recall first that by the Abhyankar lemma  (see e. g. \cite{sti}, Theorem 3.9.1)   for every point $t_0$ of $\f E_j$ the equality  
\be \l{abj} e_{R_j}(t_0)={\rm lcm} \big(e_{P}(V(t_0)), e_{W}(U(t_0))\big) \ee holds. 
In particular, $e_{R_j}(t_0)=1$ whenever $R_j(t_0)$ is not a critical value of $P$ or $W$.

If the fiber product of $P$ and $W$ consists of a unique component, then for the map $R$ defined by  \eqref{dq} we have 
\be \l{imp} {\rm c}(R)={\rm c}(P)\cup {\rm c}(W),\ee where ${\rm c}(F)$ denotes the set of critical values of 
a holomorphic map  $F$. 
Let \linebreak $\{z_1, z_2, \dots , z_r\}$ be points of the set \eqref{imp}.  
 We denote by
$(p_{i,1},p_{i,2}, ... , p_{i,u_i})$, \linebreak $1\leq i \leq r,$
the collection of local multiplicities of $P$ at the points of $P^{-1}\{z_i\}$, 
and by 
$(w_{i,1},w_{i,2}, ... , w_{i,v_i})
$, $1\leq i \leq r,$
the collection of local multiplicities of $W$ at the points of $W^{-1}\{z_i\}$. 
 The Riemann–Hurwitz formula applied to $R$  gives 
$$ \chi(\f E)=(\chi(\f C)-r)\deg P\deg W+\sum_{z\in R^{-1}\{z_1, z_2, \dots , z_r\}}1. 
$$
On the other hand, formula \eqref{abj} yields that the number of points in the preimage $R^{-1}(z_i),$ $1\leq i \leq r,$ is equal to 
$$\sum_{j_1=1}^{u_{i}} \sum_{j_2=1}^{v_{i}} \GCD(p_{i,j_1}w_{i,j_2}).$$
Thus, 
\be \l{formu} 
\chi(\f E)=(\chi(\f C)-r)\deg P\deg W+\sum_{i=1}^{r}\sum_{j_1=1}^{u_{i}} \sum_{j_2=1}^{v_{i}} \GCD(p_{i,j_1}w_{i,j_2}). 
\ee

\subsection{\l{sec} Normalizations} 
Let $F:\, \f N\rightarrow \f R$  be a holomorphic map between compact Riemann surfaces. 
Let us recall that $F$ is called a {\it Galois covering} if its automorphism group 
$$\Aut(\f N,F)=\{\sigma\in \Aut(\f N)\,:\, F\circ\sigma=F\}$$ acts transitively on  fibers  of $F$.  
Equivalently, $F$ is a Galois covering if the field extension 
$ \f M(\f N)/F^*\f M(\f R)$ is a Galois extension. In case $F$ is a Galois  covering, for the corresponding Galois group the isomorphism  
\be \l{ga} {\rm Gal}\left(\f M(\f N)/F^*\f M(\f R)\right)\cong \Aut(\f N,F)\ee holds.   
Notice that since the action of $\Aut(\f N,F)$ on fibers of $F$ has no fixed points, any element of $\Aut(\f N,F)$
is defined by its value on an arbitrary element of a fiber, 
implying that the action of $\Aut(\f N,F)$ on fibers of $F$ is transitive if and only the equality 
\be \l{dega} \vert \Aut(\f N,F)\vert =\deg F\ee 
holds. Thus,  the last equality is equivalent to the condition that $F$ is a Galois covering. 
Another equivalent condition for $F$ to be    a Galois covering is the equality \be \l{mon} \vert \Mon(F) \vert =\deg F,\ee  
where $\Mon(F)$ denotes the monodromy group of a holomorphic map $F$ (see e. g. \cite{des1}, Proposition 2.66).

Let $V:\, \f E\rightarrow \f R$  be  an arbitrary holomorphic map between compact Riemann surfaces.  Then the {\it normalization} of $V$ is defined as a compact Riemann surface  $\f N_V$ together with a holomorphic Galois covering  of the 
lowest possible degree \linebreak $\t V:\f N_V\rightarrow \f R$ such that
 $\t V=V\circ H$ for some  holomorphic map $H:\,\f N_V\rightarrow \f E$. 
The map $\t V$ is defined up to the change $\t V\rightarrow \t V \circ \alpha,$ where $\alpha\in\Aut(\f N_V)$, and is 
characterized by the property that  the field extension 
$\f M(\f N_V)/{\t V}^*\f M(\f R)$ is isomorphic to the Galois closure $\t{\f M(\f E)}/V^*\f M(\f R)$
of the extension $\f M(\f E)/V^*\f M(\f R)$. Notice that since
 $$\Mon(V)\cong {\rm Gal}\big(\t{\f M(\f E)}/V^*\f M(\f R)\big)$$ 
(see e. g. \cite{har}), 
this implies by \eqref{ga} and \eqref{dega} that the normalization of $V$ can be characterized as a Galois covering $\t V$ that factors through $V$ and satisfies the equality 
\be \l{ega} \vert \Mon(V)\vert =
\deg \t V. \ee

For a holomorphic map $V:\, \f E\rightarrow \f R$ of degree  $n\geq 2$ 
its  normalization  can be described in terms of the fiber product of 
 $V$   with itself $n$ times as follows  (see \cite{fried}, $\S$I.G, or \cite{hrs}, Section 2.2).  
For $k$, $2\leq k \leq n,$ let $\f L^{k,V}$ be an algebraic variety  consisting of $k$-tuples  of $\f E^k$ with a common image 
under $V$, 
$$\f L^{k,V}= \{(x_i)\in \f E^k\, \vert \, V(x_1)=V(x_2)=\dots =V(x_k)\},$$
and  
 $\hat{\f L}^{k,V}$  a variety obtained from $\f L^{k,V}$ by removing components that belong to the big diagonal 
$${\bf \Delta}^{k,\f E}:=\{(x_i)\in \f E^k\, \vert \, x_i=x_j \ \ {\rm for\ some} \ \ i\neq j\}$$ of  $\f E^k.$ 
Further, let  $\f L$ be  an arbitrary irreducible component of 
$\hat{\f L}^{n,V}$ and  ${\f N}\xrightarrow{\theta} \f L$ the desingularization map.  
Finally, let $\psi: {\f N}\rightarrow \f R$ be a holomorphic map induced by the composition
$$  {\f N}\xrightarrow{\theta} {\f L}\xrightarrow{\pi_i}\f E\xrightarrow{V}{\f R},$$
where  $\pi_i$ is the projection to any coordinate. In this notation, the following statement holds. 
 
\bt \l{frid} The map $\psi: {\f N}\rightarrow \f R$  is the normalization of $V$. 
\et 
\pr 
It follows  from the construction that \be \l{impl} \deg \psi=\vert \Mon(V) \vert\ee 
and that the action of $\Mon(\psi)$ on the fibers of $\psi$ has no fixed points.  The last property yields that 
$\deg \psi=\vert \Mon(\psi)\vert,$
implying that $\psi$ is a Galois covering, according to the characterization \eqref{mon}. 
Moreover, since $\psi$ factors through $V$, equality \eqref{impl} implies that $\psi$ is the normalization of $V$, according to the characterization \eqref{ega}. \qed 
\vskip 0.2cm
Notice that while the above construction is meaningless if the map $V:\, \f E\rightarrow \f R$ 
 has degree one, any such a map is a Galois covering with $\f N_V=\f E$ and $\t V=V.$

\section{ Proof of the main results}

\subsection{\l{ces} Upper and lower bounds  for $\chi(\f N_V$)} 

Let ${\f R}$ be a compact Riemann surface. We recall that an orbifold $\oo O$ on ${\f R}$ is  a ramification function $\nu:{\f R}\rightarrow \mathbb N$ which takes the value $\nu(z)=1$ except at finitely many points.  The Euler characteristic of an orbifold $\oo O=({\f R},\nu)$ is defined by the formula  
\be \l{char}  \chi(\oo O)=\chi({\f R})+\sum_{z\in {\f R}}\left(\frac{1}{\nu(z)}-1\right),\ee where $\chi({\f R})$ is the 
Euler characteristic of ${\f R}.$ 
For a holomorphic map $V:\, \f E\rightarrow \f R$ between compact Riemann surfaces, we define its {\it ramification orbifold}  
$\oo O^{ V}=(\f E,\nu)$ setting for $z\in \C\P^1$ the value $\nu(z)$     
equal to the least common multiple of local multiplicities of $V$ at the points 
of the preimage $V^{-1}\{z\}$. Notice that Theorem \ref{frid} combined with formula \eqref{abj} imply that  
\be \l{po} \oo O^{\t V}=\oo O^{ V}.\ee 

\bl \l{le1} 
Let $V:\, \f E\rightarrow \f R$ be  a holomorphic map between compact Riemann surfaces, and $\t V:\f N_V\rightarrow \f R$ its normalization. Then 
\be \l{gc} \chi(\f N_V)=\chi(\oo O^{ V})\deg \t V.\ee
\el 
\pr 
Since  $\t V:\f N_V\rightarrow \f  R$ is a Galois covering, the equality $\vert \Aut(\f N_V,\t V)\vert =\deg \t V$ holds, and   
$\t V$ is the quotient map $$\t V:\f N_V\rightarrow \f N_V/\Aut(\f N_V,\t V).$$ Thus, for any critical value  $z_i$, $1\leq i \leq r,$  of $\t V$ there exists a number $d_i$ such that $\t V^{-1}\{z_i\}$ consists of  $\deg \t V /d_i$ points, at each of which the local multiplicity of $\t V$ equals $d_i$. Applying now 
the Riemann-Hurwitz formula to $\t V$, we see that 
\be
\begin{split} 
\chi(\f N_V)&=\chi(\f R)\deg \t V -\sum_{i=1}^r\frac{\deg \t V}{d_i}\left(d_i-1\right)=\\
&=\left(\chi(\f R)+\sum_{i=1}^r\left(\frac{1}{d_i}-1\right)    \right)\deg \t V= \chi(\oo O^{\t V})\deg \t V.
\end{split}
\ee
Therefore, \eqref{gc} holds by  \eqref{po}. \qed

\bl \l{le2} Let $V:\, \f E\rightarrow \f R$ be  a holomorphic map between compact Riemann surfaces.
 Then 
\be  \l{cha} \chi(\oo O^V)\geq  \chi(\f E)+\chi(\f R)(1-\deg V).\ee
\el 
\pr It follows from the definition of $\oo O^V$  that 
\be \l{a} \chi(\oo O^V)\geq  \chi(\f R)- \vert {\rm c}(V)\vert .\ee 
On the other hand, it is clear that the number of critical values  of $V$ does not exceed the number of critical points of $V$, which in turn does not exceed the  number $\sum_{z\in \f E}(e_V(z)-1).$
Therefore, since $$\chi(\f E)=\chi(\f R)\deg V-\sum_{z\in \f E}(e_V(z)-1), $$
 we have:   
\be \l{b} \vert {\rm c}(V)\vert \leq \sum_{z\in \f E}(e_V(z)-1) =
\chi(\f R)\deg V-\chi(\f E).\ee Now \eqref{cha} follows from \eqref{a} and \eqref{b}. \qed  

 The above lemmas combined with the Hurwitz automorphisms theorem imply the following statement. 

\bt \l{mt5} 
Let $V:\, \f E\rightarrow \f R$ be  a holomorphic map between compact Riemann surfaces, and $\t V:\f N_V\rightarrow \f R$ its normalization. Then 
\be \l{ine1} \chi(\f N_V)\geq  \big(\chi(\f E)+\chi(\f R)(1-\deg V)\big)\deg \t V.\ee Furthermore, if $g(\f N_V)>1$, then 
\be \l{ine2} \chi(\f N_V)\leq -\frac{\deg \widetilde V}{42}.\ee 
\et 
\pr The first part of the theorem follows from Lemma \ref{le1} and Lemma \ref{le2}. To prove the second, we  recall that by the  Hurwitz  theorem  for a compact Riemann surface of genus $g>1$ the order of its automorphism group does not exceed $84(g-1).$  Thus, if $g(\f N_V)> 1$, then 
\be \l{ine3} 42\chi(\f N_V)\leq -\vert \Aut(\f N_V)\vert.\ee
On the other hand, since $$\Aut(\f N_V,\t V)\subseteq \Aut(\f N_V)$$ and 
the map $\t V$ is  a Galois covering, 
 it follows from \eqref{dega} that
\be \l{ine4} \deg \t V=\vert \Aut(\f N_V,\t V)\vert\leq  \vert \Aut(\f N_V)\vert. \ee
Now \eqref{ine2} follows from \eqref{ine3} and \eqref{ine4}. \qed

\subsection{\l{lii} Lifting lemma} 
Let $W:\f T\rightarrow \f C$ be a holomorphic map between compact Riemann surfaces, and 
$\f D$  a component of the  fiber product  
of the map $W$ with itself $k$, $2\leq k \leq \deg W,$ times. Then  $\f D$ is the desingularization of an irreducible  component  
$D$ of the variety ${\f L}^{k,V}$,  and  
abusing notation we will say that $\f D$ does not belong to the big diagonal of $\f T^k$, if $D$  belongs to $\hat{\f L}^{k,V}$.

Our proof of Theorem \ref{t1} and Theorem \ref{t2} is based on the following lemma of independent interest. 

\bl \l{lift}
Let 
\be \l{d1} 
\begin{CD}
\f E @>U>> \f T\\
@VV {V} V @VV W V\\ 
\f R @>P >> \f C
\end{CD}
\ee
be a reduced commutative diagram of holomorphic maps between compact Riemann surfaces such that 
$\deg V> 1$. 
 Then one can complete it to a diagram of holomorphic maps between compact Riemann surfaces
\be \l{d2} 
\begin{CD}
\f N_V @>L >> \f D\\
@VV Q V @VV F V\\ 
\f E @>U>> \f T\\
@VV {V} V @VV W V\\ 
\f R @>P >> \f C,
\end{CD}
\ee
where $\f D$ is some irreducible component of the  fiber product  
of $W:\f T\rightarrow \f C$ with itself $\deg V$ times that does not belong to the big diagonal of $\f T^{\deg V}$.  Moreover,   $V\circ Q= \t V$ and  
\be \l{lass} \deg F\leq 
(\deg W-1)\dots (\deg W-\deg V+1).\ee
Finally, if $\deg V=\deg W$,   
then $\f D=\f N_W$ and $W\circ F= \t W$.
\el 
\pr Let us set $k=\deg V$, and define the maps $$V_k: \f E^k\rightarrow \f R^k, \ \ U_k:\f E^k\rightarrow \f T^k, \ \ W_k:\f T^k\rightarrow \f C^k, \ \ P_k:\f R^k\rightarrow \f C^k$$ 
by the formulas 
\be
\begin{split} 
V_k &:\, (z_1,z_2,\dots, z_k)\rightarrow (V(z_1),V(z_2),\dots ,V(z_k)), \\
U_k &:\, (z_1,z_2,\dots, z_k)\rightarrow (U(z_1),U(z_2),\dots ,U(z_k)),\\
 W_k &:\, (z_1,z_2,\dots, z_k)\rightarrow (W(z_1),W(z_2),\dots ,W(z_k)), \\ 
P_k &:\, (z_1,z_2,\dots, z_k)\rightarrow (P(z_1),P(z_2),\dots ,P(z_k)).
\end{split}
\ee
Clearly,  the diagram 
\be \l{kaba} 
\begin{CD}
\f E^k @>U_k>> \f T^k\\
@VV {V_k} V @VV W_k V\\ 
\f R^k @>P_k >> \f C^k
\end{CD}
\ee
commutes. Furthermore, it follows from Theorem \ref{frid} that  for an arbitrary  irreducible component  $\f L$  of $\hat{\f L}^{k,V}$  the map from  $\f L$ to $U_k(\f L)$ induced  by the map 
$U_k:\f E^k\rightarrow \f T^k$  can be lifted to a map $L: \f N_V \rightarrow \f D$, where $\f D$ is the disingularization of $U_k(\f L)$, and diagram \eqref{d2} commutes for 
$$Q=\pi_{1,i}\circ \theta_{1}, \ \ \ \ F=\pi_{2,i}\circ \theta_{2},$$
where $\theta_{1}:\f N_V \rightarrow \f L$ and $\theta_{2}:\f D\rightarrow U_k(\f L)$ are the desingularization maps, and
$\pi_{1,i}:\f L\rightarrow \f E$ and $\pi_{2,i}:U_k(\f L) \rightarrow \f T$ are the projections to any coordinate.

 To prove that $\f D$  does not belong to the big diagonal of $\f T^{k}$,  it is enough to show that \be \l{mor} U_k(\hat{\f L}^{k,V})\subseteq \hat{\f L}^{k,W}.\ee  
In the notation of  Section \ref{sec},  $$\f L^{k,V}=V_k^{-1}(\Delta^{k,\f R}),$$ where $ \Delta^{k,\f R}$ is the usual diagonal in $\f R^k$, 
$$ \Delta^{k,\f R}:=\{(x_i)\in \f R^k\, \vert \, x_1=x_2=\dots =x_k\}.$$  
Since  $P_k(\Delta^{k,\f R})=\Delta^{k,\f C}$, it follows from the commutativity of \eqref{kaba} that  $$U_k({\f L}^{k,V})\subseteq {\f L}^{k,W}.$$ 

Further, 
it follows from 
\be \l{prim} V^*\f M(\f R)\cdot U^*\f M(\f T)=\f M(\f E)\ee  by the primitive element theorem that \be \l{gop} \f M(\f E)=V^*\f M(\f R)[h]\ee for some $h\in U^*\f M(\f T).$ 
As elements of  $\f M(\f E)$ separate 
points of $\f E$, 
equality \eqref{gop} implies that for all but finitely many    
$z_0\in \f R$  the map $h$ takes $\deg V$ distinct values on the set $V^{-1}\{z_0\}$.  
Since $h\in U^*\f M(\f T),$ this implies in turn that for all but finitely many    
$z_0\in \f R$
the map $U$ takes $\deg V$ distinct values on  $V^{-1}\{z_0\}$. Therefore, 
 \eqref{mor} holds. 

To finish the proof, let us observe that the degree of $F$ equals  the degree of the projection $\pi_{2,i}:U_k(\f L) \rightarrow \f T$. By \eqref{mor}, the last degree does not exceed the degree of the projection $\pi_{2,i}: \hat{\f L}^{k,W}\rightarrow \f T$, which is  equal to $(\deg W-1)\dots (\deg W-k+1)$.  
Finally, it is easy to see  that  $V\circ Q= \t V$, and
 $$\f D=\f N_W, \quad W\circ F= \t W, $$ if $\deg V=\deg W$. \qed 

\subsection{Bounds for fiber products with one component}{\it Proof of Theorem \ref{t1}.} 
If   $(\f R,P)\times_{\f C} (\f T,W)$ consists of a unique component $\f E$, then  there exists a reduced diagram \eqref{dubi} such that \be \l{fs} \deg V=\deg W.\ee
Assume first that $\deg V>1$. By Lemma \ref{lift}, we can complete diagram \eqref{dubi} to diagram  \eqref{medved},  
where $V\circ Q=\t V$ and $W\circ F= \t W$.  By the Riemann-Hurwitz formula, we have 
$$\chi(\f N_V)\leq \chi(\f N_W)\deg L,$$
implying by Theorem \ref{mt5}  that  
$$  \big(\chi(\f E)+\chi(\f R)(1-\deg V)\big)\deg \t V \leq -\frac{\deg \widetilde W\,\deg L}{42}.$$ 
Since 
\be \l{sin} \deg \t V\,\deg P= \deg  V\,\deg Q\, \deg P=\deg F\,\deg  W\,\deg L\ee
 and $W\circ F= \t W$, this implies the inequality 
\be \l{svko} \chi(\f E)+\chi(\f R)(1-\deg V) \leq -\frac{\deg P}{42},\ee 
 which  is equivalent to \eqref{svk4} by \eqref{fs}.

Assume now that $\deg V=1.$ Then $\f E\cong \f R$, $\f T\cong \f C$, and inequality \eqref{svk4} reduces to the inequality $$g(\f  R)\geq 1+\frac{\deg P}{84}.$$ On the other hand,  in case $\deg W=1$ the condition $g(\f N_W)>1$  is equivalent to the condition $\chi(\f C)\leq -2$,  and by the Riemann-Hurwitz formula we have 
 $$\chi(\f R)\leq \chi(\f C)\deg P\leq -2\deg P,$$ 
whence $$g(\f R)\geq 1+\deg P.$$ Thus, the theorem is still true although not with the best bound  
due to the fact that for $\deg V=1$ the bound \eqref{ine2} is worse than the bound $\chi(\f N_V)\leq -2$.  \qed

\br In case $\f C=\C\P^1$, Theorem \ref{t1} was also proved in the paper \cite{aol}  by a modification of the method of \cite{cur} (see \cite{aol}, Theorem 3.1). Unfortunately, by the mistake of the author, the formulation of the corresponding result in \cite{aol} was partly copied from an earlier version of the paper. As a result, it is stated in \cite{aol} that $P$ is a {\it rational function} but what is really meant is that
 $P:R\rightarrow \C\P^1$ is  a {\it holomorphic map} from  a compact Riemann surface $R$ while $W:T\rightarrow \C\P^1$ is a holomorphic map from another compact Riemann surface $ T.$
\er

\subsection{Bounds for fiber products with several components}
\noindent{\it Proof of \linebreak Theo\-rem \ref{t2}.} 
By Lemma \ref{lift},  
 we can complete diagram \eqref{dubi} to diagram  \eqref{d2} and   
arguing as in the proof of Theorem \ref{t1} we see that  
$$\chi(\f N_V)\leq \chi(\f D)\deg L$$ and   
$$  (\chi(\f E)+\chi(\f R)(1-\deg V)) \leq \frac{\chi(\f D)\,\deg L}{\deg \t V}=\frac{\chi(\f D)\,\deg P}{\deg F\,\deg W}.$$
Moreover, since the conditions of the theorem imply that $\chi(\f D)\leq -2$, we have 
\be \l{bbbb} \chi(\f E)+\chi(\f R)(1-\deg V)\leq \frac{-2\,\deg P}{\deg F\,\deg  W}.\ee
Combining now \eqref{bbbb} and \eqref{lass}, we conclude that 
$$\chi(\f E)+\chi(\f R)(1-\deg V)\leq \frac{-2\,\deg P}{\deg W(\deg W-1)\dots (\deg W-\deg V+1)}.\eqno{\Box}$$

\vskip 0.2cm
\noindent{\it Proof of Theorem \ref{rat}.}
The theorem follows from Theorem \ref{t2} for 
$W=A,$ $P=B$ and $\f R=\f T=\f C=\C\P^1,$ taking into account that for $\f R=\C\P^1$ the inequality in \eqref{svk4} is strict. Indeed,   
any holomorphic map $V:\f E\rightarrow \C\P^1$ of degree greater than one has critical values. This yields that the inequality in \eqref{a} is strict, implying that the inequalities  in \eqref{cha}, \eqref{ine1}, \eqref{svk4}, and \eqref{bbbb} are also strict.

\begin{remark}
Note that since for rational $A$ and $B$ inequality \eqref{svk4} is strict, 
for irreducible curves $E_{A,B}$ with rational $A$ and $B$, Theorem \ref{t1} gives \eqref{ma}.
\end{remark}

\vskip 0.2cm
\noindent{\it Proof of Theorem \ref{ratt}.}
Let us observe that  any irreducible
 component  $\f L$ of  ${\f L}^{k,V}$, $2\leq k\leq n,$ projects to an irreducible component $\widetilde{\f L}$ of  ${\f L}^{2,V}$, and  $g(\f L)\geq g(\widetilde{\f L})$. Moreover, if $\f L\subset \hat {\f L}^{k,V}$, then 
$\widetilde{\f L}\subset \hat{\f L}^{2,V}$. Therefore, if 
$A$ is tame, then for any $k,$ $2\leq k\leq n,$  curve 
\eqref{res}  has no component of genus zero or one that does not belong to the big diagonal of $(\C\P^1)^k$. Moreover, for any $k$, $2\leq k\leq n,$  the inequalities   
$$-k\geq -n, \ \ \ \ \ \ \ \ \ \ \frac{m}{n(n-1)\dots (n-k+1)}\geq \frac{m}{n!} 
$$
hold. Thus, if $k> 1$, then \eqref{svk11} follows from \eqref{svk10}. On the other hand, it is easy to see that if $k=1$, then $B=A\circ S$ for some rational function $S$, and  $C$ is the graph $x-S(y)=0.$ \qed

\subsection{\l{exam} On the sharpness of the bound of Theorem \ref{t1}}
In general, the bound provided by  Theo\-rem \ref{t1} is the best possible and 
is attained whenever $P$ and $W$ are the quotient  maps associated with the action of the full automorphism groups of two Hurwitz surfaces, assuming that these maps have the same ramification orbifold and the fiber product of $P$ and $W$ consists of a unique component.  Here, by a Hurwitz surface we mean a compact Riemann surface $\f R$ of genus $g>1$ with $\vert \Aut(R)\vert =84(g-1).$ 

Indeed, in terms of coverings, Hurwitz surfaces can be described as surfaces $\f R$ such that there exists a holomorphic Galois covering $P:\f R\rightarrow \C\P^1$ with the ramification orbifold $(\f O^P,\nu)$ whose 
signature is $(2,3,7)$ (the signature for which  the Euler characteristic takes the maximum possible value $-\frac{1}{42}$ among all orbifolds of negative Euler characteristic). 
By formula \eqref{gc}, for such a surface $\f R$, the genus $g(\f R)$ and the degree $n_P$ of $P$ are related by the formula 
\be \l{yaf} 84(g(\f R)-1)=n_P, \ee
and it clear that the number of points in the preimages $P^{-1}(0),$ $P^{-1}(1),$ $P^{-1}(\infty)$ is 
$$n_{1,P}=\frac{n_P}{2}, \ \ \ n_{2,P}=\frac{n_P}{3}, \ \ \ n_{3,P}=\frac{n_P}{7}$$ correspondingly.

Assuming that the fiber product of two such coverings $P$ and $W$ consists of a unique component $\f E$ and that $P$ and $W$ have the same ramification orbifold, say, defined by the equalities  \be \l{orbi} \nu(0)=2, \ \ \ \nu(1)=3,  \ \ \ \nu(\infty)=7,\ee  we obtain by formula 
\eqref{formu} that 
\begin{align} 
\chi(\f E) &=-n_P\cdot n_W+2\cdot n_{1,P}\cdot n_{1,W} +3\cdot n_{2,P}\cdot n_{2,W} +7\cdot n_{3,P}\cdot n_{3,W}=\\
&=-n_P\cdot n_W+n_P\cdot n_{1,W} +n_P\cdot n_{2,W} +n_P\cdot n_{3,W}=\\ &=-n_P\left(n_W-\frac{n_W}{2}-\frac{n_W}{3}-\frac{n_W}{7}\right)=-\frac{n_Pn_W}{42},\end{align}  whence 
$$g(\f E)=1+\frac{n_Pn_W}{84}.$$ 
On the other hand, formula \eqref{svk4} gives  
\be \l{rp} g(\f E)\geq (g(\f R)-1)( n_W-1)+1+\frac{n_P}{84},\ee 
and using  \eqref{yaf} we see that the right part of \eqref{rp} equals   
$$\frac{n_P}{84}( n_W-1)+1+\frac{n_P}{84}=1+\frac{n_Pn_W}{84}.$$ 
Thus, the equality in \eqref{svk4} is indeed attained.

In this paper, we do not make an attempt to consider the irreducibility problem for fiber products of maps of the above form in general, 
limiting ourselves  by giving an example of the irreducibility and an example of the reducibility, basing on the list of all possible genera up to  $11905$ of  Hurwitz surfaces obtained by Conder \cite{con}. The beginning of this list is 
\be \l{lis1} 
g= 3, 7, 14, 17, 118, 129, 146, 385,
411, 474, 687, 769,1009, 1025,  1459,\ee $$1537, 2091, 2131, 2185, 2663, 3404\, \dots ,
$$ 
and formula \eqref{yaf} gives 
the list of the degrees $n$ of the corresponding coverings. To save space, we present the quotients $n/84$ of these degrees
\be \l{lis2} \frac{n}{84}=g-1=2, 6, 13, 16, 117, 128, 145, 384, 410, 473, 686, 768, 1008, 1024, 1458,\ee $$ 1536, 2090, 2130, 2184, 2662, 3403\, \cdots .$$

Let us show that the fiber product of any coverings $P:\f R\rightarrow \C\P^1$ and \linebreak $W:\f T\rightarrow \C\P^1$ of the above form with $g(\f R)=3$ and $g(\f T)=14$ consists of a unique component. Indeed,  since the compositum $LK/k$ of two Galois extensions $L/k$ and $K/k$ is a Galois extension, it follows from formula \eqref{abj} that any component of the fiber product of $P$ and $W$ is a Hurwitz surface. Thus, if this fiber product  contains more than one component, then there should exist a  Hurwitz surface such that the degree of the corresponding covering is divisible by the number 
$${\rm lcm}(n^P,n^W)={\rm lcm}(84\cdot 2,84 \cdot 13)=84\cdot 2 \cdot 13,$$ but is strictly less 
than the number $n_Pn_W= 84^2\cdot 2\cdot 13$. In this case, 
list \eqref{lis2} should contain an entry that is divisible by $2\cdot 13=26$, but is strictly less than 
$84\cdot 2\cdot 13=2184.$  However, this is not true as a direct calculation shows.  

On the other hand, to show that the fiber product of any coverings $P:\f R\rightarrow \C\P^1$ and $W:\f T\rightarrow \C\P^1$ of the above form with $g(\f R)=3$ and $g(\f T)=17$ is reducible, it is enough to observe that list \eqref{lis2} does not contain an entry $$\frac{n_Pn_W}{84}=84\cdot 2\cdot 16=2688.$$

\end{document}